\documentclass {amsart}
\newtheorem{theorem}{Theorem}[section]
\newtheorem{lemma}[theorem]{Lemma}
\newtheorem{proposition}[theorem]{Proposition}

\newtheorem{assumption}[theorem]{Assumption}
\theoremstyle{definition}

\theoremstyle{remark}
\newtheorem{remark}[theorem]{Remark}

\numberwithin{equation}{section}

\usepackage{amssymb}
\usepackage{amscd}
\usepackage{amsmath}
\usepackage{graphicx}

\usepackage[all]{xy}

\newcommand{\C}{\mathbb C}
\newcommand{\ff}{\mathbb F}

\newcommand{\pp}{\mathbb P}
\newcommand{\calO}{\mathcal{O}}
\newcommand{\calE}{\mathcal{E}}

\newcommand{\st}{\overline}

\begin{document}

\title[Surfaces with $p_g=q=0$]
{Surfaces of general type with $p_g=q=0$ having a pencil of
 hyperelliptic curves of genus $3$}
\author{Giuseppe Borrelli}

\email{borrelli@mat.uniroma3.it \newline
{\em \tiny .}\ \ \ \ \ \ \ \ \ \ \ \ \ \ \ \ borrelli@dmat.ufpe.br}

\thanks{This work was supported by bolsa DTI - Instituto do Milenio/CNPq.}


\date{}

\begin{abstract}
We prove that the bicanonical map of a surfaces of general type $S$ with $p_g=q=0$ is 
non birational if there exists a pencil $|F|$ on $S$ whose general member is an hyperelliptic curve of genus $3$.
\end{abstract}

\maketitle

Let $S$ be a minimal surface of general type and let $f_g:S\dashrightarrow B$ be a rational map onto a smooth curve such
that the normalization of the general fibre is an hyperelliptic curve of genus $g$. Then the hyperelliptic involution of
the general
fibre induces a (biregular) involution $\sigma$ on $S$ and one has a map of degree two $\rho:S\dashrightarrow \Sigma$
onto a smooth ruled surface with ruling induced by $f_g$.

It is known (e.g \cite{bm},\cite{ci}) that if $g=2$ than the bicanonical map $\varphi_{2K}$ of $S$ factors through $\rho$ and moreover, if
$K_S^2\ge10$ then the non birationality of $\varphi_{2K}$ implies the existence of $f_g$ with $g=2$ (\cite{R}).
 On the other hand, we proved in
\cite{B1} that if $\varphi_{2K}$ factors through a rational map $\rho$ (generically) of degree two onto a ruled surface, 
then there exists a map $f_g$ where $g\le 4$ and, more precisely, if $g\neq 2$ then 
we have $g=3$ unless $K_S$ is ample and $q(S)=0,p_g(S)=\frac12(d-3)d+1,$ $K_S^2=2(d-3)^2$, $d=4,5$.
Finally, we recall that if $p_g(S)=0$ and $K_S^2\ge 3$, then there not exists $f_g$ with $g=2$ (\cite {X1}).

In this note we prove that if $p_g=0$ then the existence of $f_g$ with $g= 3$ implies $\varphi_{2K}$ non
birational. In particular, it follows that if $p_g=0$ and $K_S^2\ge 3$ then $\varphi_{2K}$ 
factors through a map of degree two onto a ruled surface if and only if there exists $f_g$ with $g=3$.

To motivate this work we notice that, as far as we know, all the examples of surfaces of general type with $p_g=0$ and 
non
birational bicanonical map have an $f_g$ with $g=3$, except one case when $K_S^2=3$ and $S$ is a double cover of an
 Enriques surface.

$\mathbf{Notation\ and\ conventions.}$ We work over the complex numbers. We denote by $K_S$ a canonical divisor,
by $p_g=h^0(S,\calO_S(K_S))=0$ the geometric genus and by $q=h^1(S,\calO_S(K_S))=0$ the irregularity of a
smooth (projective algebraic) surface $S$.
The symbol $\equiv$ (resp. $\sim$) will denote the linear (resp. numerical)
equivalence of divisors.
A curve on a surface has an $[r,r]$-point at $p$ if it has a point of multiplicity $r$ at $p$ which resolves
to a point of multiplicity $r$ after one blow up.

\section{surfaces with a pencil of curves of genus 3}

\begin{assumption}
Throughout the end we assume that
\begin{itemize}
\item[$a)$] $S$ is a minimal surface of general type with $p_g=q=0$ and
\item[$b)$] $f:S\dashrightarrow \pp^1$
is a rational map with connected fibres such that the normalization of general fibre $F$ is a curve of genus $g=3$.
\end{itemize}
\end{assumption}

\begin{proposition}\label{relations}
Let $F$ be a general member of $|F|$. Then,
\begin{itemize}
\item[-] either $F$ has a double point and $F\sim 2K_S$, in this case $K_S^2=1$;
\item [-] or $F$ is smooth and $F^2=2$, in this case $K_S^2=1$;
\item [-] or $F$ is smooth and $F^2\le 1$.
\end{itemize}
\end{proposition}
\begin{proof}
Suppose that $F$ has multiplicity $m_i\ge 2$ at $x_i$ and let $\st S\to S$ be the blow up of $S$ at the $x_i$'s. 
Denote by $\st F\subset \st S$ the strict transform of $F$. Then, by the adjunction formula we have 
\[
4-\st F^2-\sum m_i=2g-2-\st F^2-\sum m_i=K_{\st S}.\st F-\sum m_i=K_S.F\ge 1
\]
since $K_S$ is nef and $h^0(F)>0$. Hence, the Hodge Index Theorem says that
\[
4\ge (4-\st F^2-\sum m_i)^2=(K_S.F)^2 \ge K_S^2.F^2\ge F^2\ge\sum m_i^2\ge 4
\]
which implies that $F^2=4$ and $F$ is numerically equivalent to $2K_S$. In particular, $F$ has exactly a double point.
 
Analogously, if $F$ is smooth we get that $F^2\ge 2$ implies $F^2=2$ and $K_S^2\le 2$ with the
equality only if $F$ is numerically equivalent to $K_S$.
In the latter case $F-K_S$ is a non zero torsion element $\mu \in Pic(S)$, and 
$h^0(F)-h^1(F)+h^2(F)=\chi(F)=\chi(K_S)=1$
yields $h^1(-\mu)>0$. Therefore, the unramified covering induced by $\mu$ is irregular. A contradiction,
indeed $\pi_1^{alg}(S)\le 9$ (cfr. \cite{Re}).
\end{proof}

\section{the hyperelliptic case}
\begin{assumption}
From now on we assume that $($the normalization of$)$ the general member $F\in|F|$ is hyperelliptic.
\end{assumption}
$\mathbf{2.1.\ Involutions\ and\ double\ covers.}$
Let $\sigma_F$ be the hyperelliptic involution   
on $F$ and denote by $\sigma$ the involution induced on $S$. Then, $\sigma$ is biregular since 
$S$ is minimal of general type, and the fixed locus $Fix(\sigma)$ is union of a smooth curve $R$
and finitely many isolated points $p_1,\dots ,p_\nu$. 
\begin{lemma}\label{wei}
The base points of $|F|$ belong to the fixed locus of $\sigma$; moreover if $F^2<4$ 
$($i.e the general $F\in|F|$ is smooth$)$,
then they are distinct and isolated.
\end{lemma}
\begin{proof}
Consider $\sigma(q_1)$ for a base point $q_1$ of $|F|$.
If $F,F'$ are general members of $|F|$ we have 
\[
\sigma_F(q_1)=\sigma(q_1)=\sigma_{F'}(q_1)
\]
and so $\sigma(q_1)\in F\cap F'$.
Therefore, $q_1=\sigma(q_1)$ if $F^2=1$ or $4$.

If $F^2=2$ denote by $q_2$ the other base point ($q_2=q_1$ if $F$ and $F'$ are tangent).
 Then $\sigma(q_1)\in\{q_1,q_2\}$ and the exact sequence
\[
0\to H^0(S,\calO_S)\to H^0(S,\calO_S(F))\to H^0(F,\calO_F(q_1+q_2)) \to 0
\]
implies $h^0(F,\calO_F(q_1+q_2))=1$. Therefore, $q_1=\sigma_F(q_1)\neq q_2$ and $q_1,q_2$ are fixed points of $\sigma$.
In particular, $F$ and $F'$ meet transversally at the base points.

Suppose that there exists a component $R_{q_i}$ of $R$
passing through $q_i,i\in\{1,2\},$ so that locally $\sigma$ is the involution
$(x,y)\mapsto (x,-y)$. Notice now that since $F,F'$ meet transversally at $q_i$, they
also meet transversally $R_{q_i}$ at ${q_i}$.
A contradiction, since $F,F'$ are $\sigma$-invariant.
\end{proof}

Let $\hat S\to S$ be the blow up of $S$ at $p_1,\dots,p_\nu$ and 
denote by $\hat \sigma$ the induced involution, which is biregular since we are blowing up isolated fixed points. 
Let $\ \rho: \hat S\to \Sigma=\hat S/\hat \sigma$ be the projection onto the quotient. Hence, $\Sigma$ is a 
smooth rational surface and $\rho$ is a (finite) double cover.
Denote by $\hat B=\rho (\hat R)=\rho(\pi^*(R)+\sum E_i)$ the branch
curve of $\rho$, where $E_i=\pi^{-1}(p_i),\ i=1,\dots,\nu$. Then, $\hat B$ is a smooth curve linearly equivalent
to $2\hat \Delta$ for some $\hat \Delta\in Pic(\hat \Sigma)$.
By [3, Proposition 1.2] we have the following equalities:
\begin{align}
\nu=&\ K_S.R+4=K_S^2+4-2h^0(\hat\Sigma,2K_{\hat \Sigma}+\hat\Delta) \label{nu} \\
K_{\hat S}^2=&\ K_S^2-\nu                                       \label{K2}
\end{align}
and $\varphi_{2K_S}$ factors trough $\rho$ if and only if $h^0(\hat\Sigma,2K_{\hat \Sigma}+\hat\Delta)=0$.
Notice that $\nu\ge 4$ since $K_S$ is $nef$.
\begin{lemma}                                                                            \label{K1}
Assume $K_S^2=1$, then $\varphi_{2K_S}$ factors trough $\rho$.
\end{lemma}
\begin{proof}
Otherwise it would be $h^0(\hat\Sigma,2K_{\hat \Sigma}+\hat\Delta)>0$ and hence $\nu\le 3$.
\end{proof}
From now on, in this section we assume $K_S^2\ge 2$.
Therefore, $F$ is smooth and, by Lemma \ref{wei}, $\Sigma$ is birationally ruled by $|\hat \Gamma|$, where
$\hat \Gamma$ is the image of $F$.
Let $\omega :\Sigma \to \ff_e,\ e\ge 0$ be a birational morphism 
 such that $\Gamma=\omega_\ast(\hat\Gamma)$ is a ruling, and
consider a factorization
$\omega=\omega_d\circ\dots\circ\omega_1$ where $\omega_i:\Sigma_{i-1}\to\Sigma_i$ is the blow up at
$q_i\in\Sigma_i, i\ge 1, \Sigma_0:=\Sigma$ and $\Sigma_d=\ff_e$.

Denote by $\calE_i,\calE_i^\ast$ respectively the
exceptional curve of $\omega_i$ and its total transform on $\Sigma$ and  by
$B=\omega_\ast(\hat B)$ denote the image of
$\hat B$ on $\ff_e$.
Since $\hat B\equiv 2\hat\Delta$ we have $ B\equiv 2\Delta$, where $\Delta=\omega_\ast(\hat \Delta)\in Pic(\ff_e)$
and hence $B\equiv 8C_0+2b\Gamma$, where $C_0$ is the $(-e)$-section. Recall that
 $K_{\ff_e}\equiv -2C_0 -(2+e)\Gamma$ and $K_\Sigma\equiv \omega^\ast(K_{\ff_e})+\sum\calE_i^\ast$.

Finally, notice that the branch curve $\hat B$ contains exactly $\nu$ $(-2)$-curves (they correspond to the isolated fixed points
of $\sigma$). In particular, the $(-2)$-curves arising from the base points of $F$ map to sections of 
$\Sigma_e\to \pp_1$, while each one of the others maps either to a point or to a fibre $\Gamma$.

\begin{lemma}\label{known}
In the above situation
\begin{itemize}
\item[$i$)] $\hat S$ is the canonical resolution of the double cover of $\ff_e$ branched along $B$;
\item[$ii$)] we can assume that:
\begin{itemize}
\item [$ii.a$)] the essential singularities of $B$ are at most $[5,5]$-points, and
\item [$ii.b$)]if $[x'\to x]$ is a $[5,5]$-point then the fibre $\Gamma_x$ through $x$ belongs to $B$;
\end{itemize}
\item[$iii)$] 
 assume $(ii)$, then
\begin{itemize}
\item [$iii.a$)] 
there is a $1$ to $1$ correspondence 
\[
\left\{
\begin{array}{c}
(-2)\textrm{-curves} \ contained\\ in\ \hat B\ contracted\ to\ points
\end{array}
\right\}
\longleftrightarrow
\left\{
\begin{array}{c}
[r,r]\textrm{-points} \\
of\ B
\end{array}
\right\}
\]
and
\item [$iii.b)$] 
there is a $1$ to $1$ correspondence 
\[
\left\{
\begin{array}{c}
(-2)\textrm{-curves} \ contained\\ in\ \hat B\ which\ map\ to\ fibres
\end{array}
\right\}
\longleftrightarrow
\left\{
\begin{array}{c}
fibres\ passing\ through\\ a\ [5,5]\textrm{-point of } B
\end{array}
\right\}
\]
\end{itemize}

\end{itemize}
\end{lemma}
\begin{proof}
For $i)$ see [3, Lemma 1.3], while for $ii)$ and $iii)$ see [9, Lemmas 6 and 7].
\end{proof}

\begin{proposition}\label{relaz}
Assume that $\omega$ has the property $ii)$ of Lemma \ref{known}. Denote by $a_1$ (resp. $a_2,a_3$) the number of
$[5,5]$-points $($resp. $[3,3]$-points, $4$ and $5$-tuple points$)$ of $B$.
Then 
\begin{align*}
d\ge \ 2a_1+&\ 2a_2+a_3; \  \  \  \nu-F^2=2a_1+a_2 \\
\frac32& K_{\hat S}^2+12=a_1+a_2+a_3 
\end{align*}
\end{proposition}

\begin{proof}
Since $\hat B$ is smooth, the inequality is clear, just notice that to resolve 
a singularity of type [$r,r$] we need to blow up twice.
The equality $\nu-F^2=2a_1+a_2$ follows from $ii,b)$ and $iii)$ of the above lemma.

By Lemma \ref{known} and \cite{H}, we have $\hat B=\omega^\ast(B)-\sum 2[\frac{m_i}2]\calE_i^\ast$
where $m_i$ is the multiplicity of $B_i$ at $q_i$.
 Hence, since $\rho$ is a double cover we get
\begin{align*}
\chi(\hat S)=&\ \chi(\Sigma)+\chi(K_\Sigma+\hat\Delta)
=\frac12(K_{\ff_e}+\Delta).\Delta+2\chi(\ff_e)-\frac12\sum [\frac{m_i}2]([\frac{m_i}2]-1)
 \\
K_{\hat S}^2=&\ 2(K_\Sigma+\hat\Delta)^2=2(K_{\ff_e}+\Delta)^2-\sum([\frac{m_i}2]-1)^2
\end{align*}
and hence
\begin{align*}
1=&\ \frac12(6b-12e-8)+2-\frac12(a_18+a_22+a_32)= \\
=&\ 3b-6e-2-4a_1-a_2-a_3; \\
K_{\hat S}^2=&\ 2(4b-8e-8)-2(a_15+a_21+a_31)= \\
=&\ 8b-16e-16-10a_1-2a_2-2a_3
\end{align*}
because $B$ has at most $[5,5]$-points.
Therefore,
\[
3K_{\hat S}^2-8=3(8b-16e-16-10a_1-2a_2-2a_3)-8(3b-6e-2-4a_1-a_2-a_3)
\]
and so
\[
\frac32K_{\hat S}^2+12=a_1+a_2+a_3\\
\]
\end{proof}
\begin{remark}
Because of the injection $H^0(\Sigma,\C) \hookrightarrow H^0(\hat S,\C)$ induced by $\rho$ we have
 $K_\Sigma^2\ge K^2_{\hat S}$. Therefore, $d\le 8-K_{\Sigma}^2\le 8-K_{\hat S}^2$
since $\omega$ is a sequence of $d$ blow ups.
We set $\hat d:=8-K_{\hat S}^2$.
\end{remark}
\begin{lemma}
We have $[\frac{\nu-F^2}2]\ge a_1\ge\nu-F^2-\frac12\hat d$
\end{lemma}
\begin{proof}
The first inequality follows from Proposition \ref{relaz}. 
For the latter suppose that $a_1=\nu-F^2-\frac12\hat d-i, i\ge 1$. Then we have 
$a_2=\nu-F^2-2a_1=\hat d+2i-\nu+F^2$ and so $d\ge 2a_1+2a_2=\hat d+2i>\hat d$.
A contradiction.
\end{proof}
\begin{remark}
Notice that if $\omega$ has the property $ii)$ of Lemma \ref{known} and $B$ has a $[5,5]$-point at $p$,
then each irreducible component of $B$ passing through $p$ is tangent to the fibre which contains $p$. In particular,
if $F^2>0$ then $B$ contains $F^2$ sections of
$\ff_e\to \pp_1$ and hence it has no $[5,5]$-points, i.e $a_1=0$.
\end{remark}

\section{the main results}

\begin{theorem}\label{main}
Let $S$ be a surface of general type with $p_g(S)=q(S)=0$. Suppose that there exists a rational map 
$S\dashrightarrow \pp_1$ such that the normalization of the general fibre $F$ is an hyperelliptic curve of genus 3.
Then,
\begin{itemize}
\item[$i)$] the bicanonical map of $S$ is composed with the hyperelliptic involution.
\item[$ii)$] $K_S^2\le 8$, and $K_S^2\le 3$ if $F^2=1$. 
\end{itemize}
\end{theorem}
\begin{proof}
By Lemma \ref{K1} and Proposition \ref{relations} we can assume $F^2\le 1$.  By Proposition \ref{relaz}, we have
\[
\frac32K_{\hat S}^2+12=a_1+a_2+a_3
\]
and by the above lemma we can write
$a_1=\nu-F^2-\frac12\hat d+i,$ where 
\[ 0\le i\le \left[\frac{\nu-F^2} 2\right]-\nu+F^2+\frac12\hat d\le \frac12(\hat d-\nu+F^2)\]
Then, $a_2=\nu-F^2-2a_1=\hat d-\nu+F^2-2i$ and $a_3\le \hat d -2a_1-2a_2=2i$. Therefore, we get
\begin{align*}
3K_{\hat S}^2+24=&\ 2\left(\frac {\hat d}2-i+a_3\right)\le \hat d+2i \le \\
\le &\ 2 \hat d-\nu + F^2 =16-2K_{\hat S}^2-\nu+ F^2
\end{align*} 
and hence, 
\[
4\le \nu\le -5K_{\hat S}^2-8+F^2=12+F^2-10h^0(\Sigma,\calO_\Sigma(2K_\Sigma+\hat \Delta))\le 
13-10h^0(\Sigma,\calO_\Sigma(2K_\Sigma+\hat \Delta))
\]
by \ref{nu} and \ref{K2}. It follows that $h^0(\Sigma,\calO_\Sigma(2K_\Sigma+\hat \Delta))=0$
and whence $\varphi_{2K}$ factors through the involution.

Finally, notice that:
\begin{itemize}
\item [$a)$]
if $F^2=0$ then $K_S^2+4=\nu \le 12$, i.e. $K_S^2\le 8$;
\item[$b)$] if $F^2=1$ then $a_1=0$, hence
\[
12=2a_1+2a_2+2a_3=2(\nu-1)+2a_3=2K_S^2+6+2a_3
\]
and so
\[
K_S^2=3-a_3\le 3
\]
\end{itemize}
In particular, $K_S^2=9$ does not occur.  \end{proof}
An involution $\iota$ on $S$ is $rational$ if the quotient $S/\iota$ is a rational surface. In the situation of Section
$2$ we have the following commutative diagram
\[
 \xymatrix@!0{
 \Hat S\ar [dd]_{\rho} \ar[rrr] &&& S  \ar[dd]   \\
 \\
\Hat \Sigma \ar[rrr] & & & S/\sigma 
}
\]
where $S/\sigma$ is rational and, in particular, $\varphi_{2K}$ factors through $\sigma$ if and only if
factors through $\rho$.

\begin{theorem}
Let $S$ be a minimal surface of general type with $p_g=q=0$. Then the bicanonical map of $S$ factors through
 a rational involution
if and only if there exists a map $f_g:S\dashrightarrow \pp_1$such that the normalization of the general fibre
is an hyperelliptic curve of genus $g\le3$ $(g=3$ if $K_S^2\ge 3)$. 
\end{theorem}
\begin{proof}
Assume that there exists $f_g$. If $g=2$ see \cite{bm}, \cite{ci} (in fact, the bicanonical system cuts on the general
fibre a subserie of the $g^1_2$). If $g=3$ then we have Theorem \ref{main}.
Conversely, if the bicanonical map factors through a rational involution then by \cite{B1} there exists $f_g$ with $g\le 3$.
Finally, Xiao G. [9, Theoreme 2] 
proved that if there exists $f_g$ with $g=2$ then $K_S^2\le 2$.

\end{proof}

\vspace{10 mm}

{ \scshape Giuseppe Borrelli}

\textsc {\small Departamento de Matematica}

\textsc {\small Universidade de Pernambuco}

\textsc {\small Cidade Universitaria}

\textsc {\small 50670-901 Recife-PE}

\textsc {\small Brasil}

\end{document}